\documentclass[12pt,draftcls,onecolumn]{IEEEtran}

\usepackage{amsfonts}

\IEEEoverridecommandlockouts                              
\usepackage{amsmath}
\usepackage{amssymb}

\title{Converse Lyapunov-Krasovskii Theorems for Systems Described by Neutral Functional Differential Equation in Hale's Form }

\author{
P. Pepe {\it \ \ }
\hspace*{.2ex} \hspace*{.0ex}
\ \ \ \ \  \ I. Karafyllis {\it \ \ }
\hspace*{.2ex} 
\hspace*{.0ex} 
\thanks{
The work of P. Pepe was supported by the Center of Excellence for
Research DEWS, L'Aquila, Italy and by the Italian MIUR Project PRIN 2009.
}
\thanks{P. Pepe is with
the Department of Information Engineering, Computer Science, and Mathematics, 
University of L'Aquila, L'Aquila, Italy,
{\tt\small pierdomenico.pepe@univaq.it}.}%
\thanks{I. Karafyllis is with
the Department of Environmental Engineering, Tehcnical University of Crete, Chania, Greece,
{\tt\small iasson.karafyllis@enveng.tuc.gr}.
}%
}

\newtheorem{definition}{Definition}
\newtheorem{theorem}{Theorem}

\def\bn{\begin{eqnarray*}}
\def\be{\begin{eqnarray}}
\def\en{\end{eqnarray*}}
\def\ee{\end{eqnarray}}
\def\ba{\begin{array}}
\def\ea{\end{array}}
\def\bm{\left[\begin{array}}
\def\em{\end{array}\right]}

\begin{document}
\addtolength{\baselineskip}{-.08mm}

\maketitle
\thispagestyle{empty}   
\pagestyle{empty}       

\begin{abstract}
In this paper we show that the existence of a Lyapunov-Krasovskii functional is necessary and sufficient condition for the uniform global asymptotic stability and the global exponential stability of  time-invariant systems described by neutral functional differential equations in Hale's form. It is assumed that the difference operator is linear and strongly stable, and that the map in the right-hand side of the equation is Lipschitz on bounded sets. A link between global exponential stability and input-to-state stability is also provided. 
\noindent
{\it @ The extended version of this paper has been submitted to the International Journal of Control, Taylor \& Francis.}

\end{abstract}

\IEEEpeerreviewmaketitle

\section{Introduction}

In (Karafyllis, 2006), (Karafyllis, Pepe \& Jiang, 2008 (A)), (Karafyllis, Pepe \& Jiang, 2008 (B)), (Karafyllis \& Jiang, 2010), 
many converse Lyapunov theorems have been presented, for many global stability notions, for systems described by retarded functional differential equations, in a very general setting. Time varying delays, disturbances, time-varying equations are considered. Besides the global uniform asymptotic stability, also the input-to-state stability and the weighted input-to-output stability are investigated. The reader can refer to the recent monograph (Karafyllis \& Jiang, 2011)  
for an extensive presentation of this topic. As far as systems described by neutral equations are concerned, converse Lyapunov-Krasovskii theorems are available for the (local, uniform) asymptotic stability (see Cruz \& Hale, 1970, Ben\'a \& Godoy, 2001). Methods for constructing the quadratic Lyapunov-Krasovskii functional are also provided, for the linear case, in (Rodriguez, Kharitonov, Dion \& Dugard, 2004, Kharitonov, 2005, Kharitonov, 2008, Velazquez-Velazquez \& Kharitonov, 2009). Converse Lyapunov-Krasovskii methods are also used for establishing instability criteria for linear neutral systems in (Mondi\'e, Ochoa \& Ochoa, 2011). As far as global asymptotic stability notions are concerned, to our knowledge, converse Lyapunov-Krasovskii theorems for nonlinear neutral systems are not yet available in the literature. 

In this paper, we consider time-invariant systems described by neutral equations in Hale's form (see Hale \& Lunel, 1993, Kolmanovskii \& Myshkis, 1999). The difference operator can involve an arbitrary number of arbitrary discrete time-delays. It is assumed to be linear and strongly stable (see Hale \& Lunel, 1993). The map in the right-hand side of the equation is assumed to be Lipschitz on bounded sets. An arbitrary number of arbitrary discrete and distributed time-delays can appear in the right-hand side of the equation. 
We prove here that the well known conditions, extended to the whole state space, in the Lyapunov-Krasovskii Theorem for the (uniform) local asymptotic stability of the origin (see Kolmanovskii \& Nosov, 1982, Kolmanovskii \& Nosov, 1986, Hale \& Lunel, 1993, Kolmanovskii \& Myshkis, 1999), are not only sufficient, but also necessary. Moreover, we show here converse Lyapunov-Krasovskii theorems for global exponential stability, and a link between global exponential stability and input-to-state stability.

\medskip

\noindent

\bigskip

\noindent
{\bf Notations.} $R$ denotes the set of real numbers, $R^{\star}$ denotes
the extended real line $[-\infty,+\infty]$, $R^+$ denotes the set
of non negative reals $[0,+\infty)$, $Z^+$ denotes the set of integers in $R^+$. For $s\in R^+$, $[s]$ denotes the largest number in $Z^+$, smaller or equal to $s$. The symbol $\vert \cdot
\vert$ stands for the Euclidean norm of a real vector, or the
induced Euclidean norm of a matrix. For a positive integer $n$, for a positive real $\Delta$
(maximum involved time-delay): ${\mathcal C}$ 
denotes the Banach space of the continuous functions mapping $[-\Delta,0]$
into $R^n$, endowed with the supremum norm, indicated with the symbol $\Vert \cdot \Vert$; $W^{1,\infty}$ denotes the space of the absolutely continuous functions in ${\mathcal C}$, with essentially bounded derivative. 
With the symbol $\Vert \cdot \Vert_a$ is denoted any semi-norm in ${\mathcal C}$ (see Pepe \& Jiang, 2006, Yeganefar, Pepe \& Dambrine, 2008).
For a function $x:[-\Delta, c)\to R^n$, with $0<c\le
+\infty$, for any real $t\in [0,c)$, $x_t$ is the function in
${\mathcal C}$ defined as $x_{t}(\tau)=x(t+\tau)$, $\tau \in
[-\Delta,0]$. 
For positive real $H$, $\phi\in {\mathcal
C}$, $C_{H}(\phi)$ denotes the subset (of ${\mathcal C}$) $\{\psi \in {\mathcal C}: \Vert
\psi-\phi\Vert\le H\}$. For $C_H$, $C_H(0)$ is meant.   
For positive integer $m$, positive real $\delta$, $B_{\delta}$ denotes the subset (of $R^m$) $\{u\in R^m: \vert
u\vert\le\delta\}$.
For positive integer $n$, a map $Q: {\mathcal C}\to R^n$ is said to be: locally Lipschitz if, for any $\phi$ in ${\mathcal C}$, there exist positive reals $H$, $L$ such that, for any $\phi_1, \phi_2 \in C_H(\phi)$, the 
inequality $\vert Q(\phi_1)-Q(\phi_2)\vert \le L \Vert
\phi_1-\phi_2\Vert$ holds; Lipschitz on bounded sets if, for any positive real $H$, there exists a positive real $L$ such that, for any $\phi_1, \phi_2 \in {\mathcal C}_H$,  the 
 inequality $\vert Q(\phi_1)-Q(\phi_2)\vert \le L\Vert
 \phi_1-\phi_2\Vert$ holds;
globally Lipschitz if there exists a positive real $L$ such that, for any $\phi_1, \phi_2 \in {\mathcal C}$, the 
inequality $\vert Q(\phi_1)-Q(\phi_2)\vert \le L \Vert
\phi_1-\phi_2\Vert$ holds.
For positive integers $n$, $m$, a map $Q: {\mathcal C}\times R^m\to R^n$ is said to be Lipschitz on bounded sets if, for any positive reals $H$, $\delta$, there exists a positive real $L$ such that, for any $\phi_1, \phi_2 \in {\mathcal C}_H$, for any $u_1, u_2\in B_{\delta}$,  the 
 inequality $\vert Q(\phi_1,u_1)-Q(\phi_2,u_2)\vert \le L\left (\Vert
 \phi_1-\phi_2\Vert+\vert u_1-u_2\vert\right )$ holds.
For positive integer $m$, a Lebesgue measurable function $v:R^+\to R^m$ is said to be essentially bounded if
$ess\sup_{t\ge 0}\vert v(t) \vert<\infty$. The essential supremum norm of a Lebesgue measurable and 
essentially bounded function is indicated again with the symbol $\Vert
\cdot \Vert$.  For given times $0\le
T_1<T_2$, we indicate with $v_{[T_1,T_2)}:R^+\to R^m$ the function
given by $ v_{[T_1,T_2)}(t)= v(t) $ for all $t \in [T_1,T_2)$ and
$=0$ elsewhere. An input $v$ is said to be locally
essentially bounded if, for any $T>0$, $v_{[0,T)}$ is essentially
bounded.
 Let us
here recall that a function $\gamma:R^+\to R^+$ is: of class
${\cal P}$ if it is continuous, zero at zero, and positive for any
positive real; of class ${\cal K}$ if it is of class ${\cal P}$
and strictly increasing; of class ${\cal K}_{\infty}$ if it is of
class ${\cal K}$ and it is unbounded; of class ${\cal L}$ if it is
continuous and it monotonically decreases to zero as its argument
tends to $+\infty$. A function $\beta:R^+\times R^+\to R^+$ is of
class ${\cal KL}$ if $\beta(\cdot, t)$ is of class ${\cal K}$ for
each $t\ge 0$ and $\beta(s,\cdot)$ is of class ${\cal L}$ for each
$s\ge 0$. 

\medskip

\noindent
{\bf Acronyms.} Throughout the paper, as standard,
ODE stands for ordinary differential equation, RFDE stands for retarded functional differential equation,
NFDE stands for neutral functional differential equation, CTDE stands for continuous time difference equation, GAS stands for global asymptotic stability or
globally asymptotically stable, GES stands for global exponential stability or globally exponentially stable, ISS stands for input-to-state stability or input-to-state stable. 

\section {NFDEs and Converse Lyapunov-Krasovskii Theorems for 0-GAS and 0-GES Properties}

Let us consider the following NFDE in Hale's form (see Hale \& Lunel, 1993, Kolmanovskii \& Myshkis, 1999)

\begin{eqnarray}\label{NFDE}&&
\frac{d}{dt} {\mathcal D}x_t=f(x_t), \qquad t\ge 0,  \nonumber
\\ && x(\tau)=\xi_0(\tau),\qquad \tau \in [-\Delta,0], \qquad   \xi_0\in {\mathcal C},
\end{eqnarray}
where: $x(t)\in R^n$, $n$ is a positive integer; $\Delta>0$ is the
maximum involved time-delay; the map $f:{\mathcal C}\to R^n$ is Lipschitz on bounded sets and satisfies $f(0)=0$; the operator
$\mathcal D:{\mathcal C}\to R^n$ is defined, for $\phi \in {\mathcal C}$, as
\begin{equation}\label{defDoperator}
{\mathcal D}\phi=\phi(0)-\sum_{j=1}^pA_j\phi(-\Delta_j),
\end{equation}
with $p$ a positive integer, $\Delta_j$ positive reals satisfying $\Delta_j\le \Delta$, $j=1,2,\dots, p$, $A_j$ matrices in $R^{n\times n}$, $j=1,2,\dots,p$. 
It is assumed that the operator ${\mathcal D}$ is strongly stable (see Definition 6.2, pp. 284, in Hale \& Lunel, 1993). 
Let us here recall that: in the case $p=1$, the operator ${\mathcal D}$ is strongly stable if and only if the eigenvalues of $A_1$ are located inside the open unitary disk; in the case of multiple delays, a necessary and sufficient condition for the strong stability of the operator ${\mathcal D}$ is provided in Theorem 6.1, pp. 284-287, in (Hale \& Lunel, 1993). Sufficient conditions for the strong stability of the operator ${\mathcal D}$, in terms of Matrix Inequalities, are provided in (Pepe \& Verriest, 2003), (Pepe, 2005), (Gu \& Liu, 2009),
(Gu, 2010), (Li \& Gu, 2010), (Li, 2012).

\noindent
In the following well known definition, reported here for reader's convenience, global asymptotic stability is meant uniform with respect to ${\mathcal C}_H$, for any positive real $H$ (see Definition 4.4, pp. 149-150, in Khalil, 2000, for systems described by ODEs, see also Definition 1.1, pp. 130-131 in Hale \& Lunel, 1993, as far as the local, uniform asymptotic stability of systems described by RFDEs is concerned).
\begin{definition}\label{def0GAS}
The system described by (\ref{NFDE}) is said to be $0$-GAS if:

\begin{itemize} 
\item[i)] for any $\epsilon >0$ there exists $\delta>0$ such that, for any initial condition $\xi_0\in \mathcal C_{\delta}$, the solution exists for all $t\in R^+$ and, furthermore, satisfies
\begin{equation}\vert x(t)\vert<\epsilon, \qquad t\ge 0;
\end{equation}
moreover, $\delta$ can be chosen arbitrarily large for sufficiently large $\epsilon$;
\item[ii)] for any positive real $H$, for any positive real $\epsilon$, there exists a positive real $T$ such that, for any initial condition $\xi_0\in \mathcal C_{H}$, the corresponding solution exists for all $t\in R^+$ and, furthermore, satisfies
\begin{equation}
\vert x(t)\vert < \epsilon, \qquad \forall \ t\ge T
\end{equation}
\end{itemize}
\end{definition}

\noindent Property i) in Definition \ref{def0GAS} is equivalent to Lyapunov and Lagrange stability for dynamical systems. The following well known definition concerns the global exponential stability (see Definition 4.5, p. 150, in Khalil, 2000, as far as systems described by ODEs are concerned, see Krasovskii, 1963, as far as systems described by RFDEs are concerned).
\begin{definition}\label{definexpMlambda} The system described by (\ref{NFDE}) is said to be $0$-GES if there exist positive reals $M$, $\lambda$ such that, for any initial condition $\xi_0\in \mathcal C$, the corresponding solution of (\ref{NFDE}) exists for all $t\in R^+$ and, furthermore, satisfies the inequality
\begin{equation}\label{expwithM}
\vert x(t)\vert \le Me^{-\lambda t}\Vert \xi_0\Vert, \qquad t\ge 0
\end{equation}
\end{definition}
\noindent
Notice, in Definition \ref{definexpMlambda}, that $M\ge 1$ is mandatory.
For a locally Lipschitz functional
$V:{\mathcal C}\to R^+$, the derivative of the
functional $V$, $D^+V:{\mathcal C}\to R^{\star}$, is defined (in Driver's form,
see Driver, 1962, Pepe, 2007 (A)), for $\phi\in {\mathcal C}$, as
\begin{equation}\label{derivatatotale}
 D^+ V(\phi)=
\limsup_{h\to 0^+} {1\over h} \left ( V(\phi_{h})-V(\phi) \right
),
\end{equation}
where: for $0<h<\Delta$, $\phi_h \in {\mathcal C}$ is given
by
\begin{eqnarray}\label{fictitioussystemneutralold}  \phi_h (s) =  \left \{ \begin{array}{cc} \phi(s+h),
\qquad s\in [-\Delta,-h];     \\ {\cal
D}\phi+f(\phi)(s+h) -{\cal
D}\phi_{s+h}^{\star}+\phi(0),   \qquad  s\in (-h,0];
\end{array}
\right .
\end{eqnarray}
for $0<\theta\le h$, $\phi^{\star}_{\theta} \in {\mathcal C}$ is given by
\begin{equation}\label{phithetastarold}
\phi^{\star }_{\theta}  (s) =  \left \{ \begin{array}{cc}
\phi(s+\theta), & s\in [-\Delta,-\theta];   \\  \phi(0), &  s\in
(-\theta,0]
 \end{array} \right .
\end{equation}

\noindent
Let us here recall that, for any locally Lipschitz functional $V:{\mathcal C}\to R^+$, the following result holds (see Driver, 1962, Pepe, 2007 (A))
\begin{equation}\label{krasovskiiequivdriver}
\limsup_{h\to 0^+}\frac{V(x_{t+h})-V(x_t)}{h}=D^+V(x_t), \qquad t\in [0,b),
\end{equation}
where $x_t$ is the solution of $(\ref{NFDE})$ in a maximal time interval $[0,b)$, $0<b\le+\infty$.

\noindent
The main result of this paper is given by the necessity part of the following Theorem (see Theorem 8.1, p. 293, in Hale \& Lunel, 1993, and Theorem 7.2 in Cruz \& Hale, 1970, as far as local, uniform asymptotic stability is concerned).

\begin{theorem}\label{converseGAS} The system described by (\ref{NFDE}) is $0$-GAS if and only if there exist a locally Lipschitz functional $V:{\mathcal C}\to R^{+}$, functions $\alpha_1$, $\alpha_2$ of class $\mathcal K_{\infty}$, a function $\alpha_3$ of class ${\mathcal K}$, such that the following conditions hold for all $\phi \in {\mathcal C}$:
\begin{itemize}
\item[i)] $\alpha_1\left (\vert{\mathcal D}\phi\vert\right )\le V(\phi)\le \alpha_2\left (\Vert\phi\Vert\right )$;
\item[ii)] $D^+V(\phi)\le -\alpha_3 \left (\vert{\mathcal D}\phi\vert\right )$
\end{itemize}
\end{theorem}
\noindent
The following Theorems provide necessary and sufficient Lyapunov-Krasovskii conditions for the $0-$GES property.

\begin{theorem}\label{karconverseGES} The system described by (\ref{NFDE}) is $0$-GES if and only if there exist a locally Lipschitz functional $V:{\mathcal C}\to R^{+}$ and positive reals $a_1$, $a_2$, $a_3$, such that the following conditions hold for all $\phi \in {\mathcal C}$:
\begin{itemize}
\item[i)] $a_1\vert{\mathcal D}\phi\vert \le V(\phi)\le a_2\Vert \phi\Vert$;
\item[ii)] $D^+V(\phi)\le -a_3 V(\phi)$
\end{itemize}
\end{theorem}

\begin{theorem}\label{converseGESflip} If the map $f$ in (\ref{NFDE}) is globally Lipschitz, then the system described by (\ref{NFDE}) is $0$-GES if and only if there exist a globally Lipschitz functional $V:{\mathcal C}\to R^{+}$, 
such that the following conditions hold for all $\phi \in {\mathcal C}$:
\begin{itemize}
\item[i)] $a_1\vert{\mathcal D}\phi\vert \le V(\phi)\le a_2\Vert \phi\Vert$;
\item[ii)] $D^+V(\phi)\le -a_3 V(\phi)$
\end{itemize}

\end{theorem}

\noindent

The following
Theorem will be used for establishing a link between $0-$GES and ISS properties in next section.

\begin{theorem}\label{converseGEScor} If the map $f$ in (\ref{NFDE}) is globally Lipschitz, then the system described by (\ref{NFDE}) is $0$-GES if and only if there exist a globally Lipschitz functional $V:{\mathcal C}\to R^{+}$, a semi-norm $\Vert \cdot \Vert_a$ and positive reals $a_1$, $a_2$, $a_3$, $a_4$, such that the following conditions hold for all $\phi \in {\mathcal C}$: 
\begin{itemize}
\item[i)] $a_1\vert{\mathcal D}\phi\vert \le V(\phi)\le a_2\Vert \phi \Vert_a$;
\item[ii)] $D^+V(\phi)\le -a_3 \Vert \phi \Vert_a$;
\item[iii)]  $\Vert \phi \Vert_a\le a_4 \Vert \phi\Vert$.
\end{itemize}
\end{theorem}
\noindent

\section{A Link between 0-GAS and ISS Properties for Systems Described by NFDEs}\label{sectionlinkiss}

\noindent
A link between $0-$GES and ISS properties is here provided (see Lemma 4.6, pp. 176-177, in Khalil, 2000, as far as $0$-GES systems described by ODEs are concerned, see 
Yeganefar, Pepe \& Dambrine, 2008, as far as $0-$GES systems described by RFDEs are concerned). 
Let us consider the following NFDE in Hale's form 

\begin{eqnarray}\label{NFDEinput}&&
\frac{d}{dt} {\mathcal D}x_t=f(x_t,u(t)), \qquad t\ge 0,  \nonumber
\\ && x(\tau)=\xi_0(\tau),\qquad \tau \in [-\Delta,0], \qquad   \xi_0\in {\mathcal C},
\end{eqnarray}
where: $x(t)\in R^n$, $n$ is a positive integer; $\Delta>0$ is the
maximum involved time-delay; the map $f:{\mathcal C}\times R^m\to R^n$ is Lipschitz on bounded sets and satisfies $f(0,0)=0$; $m$ is a positive integer; $u(\cdot)$ is a Lebesgue measurable, locally essentially bounded input signal; the operator
$\mathcal D:{\mathcal C}\to R^n$, defined as in (\ref{NFDE}), is assumed to be strongly stable. 
It is assumed that there exists a map $\overline f:R^n\times {\mathcal C}\times R^m\to R^n$, independent of the second argument at $0$ (see Pepe, 2011, see Definition 5.1, p. 281, in Hale \& Lunel, 1993), such that, for any $\phi\in {\mathcal C}$, $u\in R^m$, the equality $f(\phi,u)=\overline f(\phi(0),\phi,u)$ holds.

\begin{definition} (see Sontag, 1989)
The system described by (\ref{NFDEinput}) is said to be ISS if there exists a function $\beta$ of class ${\mathcal {KL}}$ and a function $\gamma$ of class ${\mathcal {K}}$ such that, for any initial condition $\xi_0\in {\mathcal C}$, for any Lebesgue measurable, locally essentially bounded input signal $u$, the corresponding solution of (\ref{NFDEinput}) exists for all $t\ge 0$ and, furthermore, satisfies the inequality
\begin{equation}
\vert x(t)\vert\le \beta(\Vert \xi_0\Vert,t)+\gamma (\Vert u_{[0,t)} \Vert), \qquad t\in R^+
\end{equation}
\end{definition}

\noindent
For a locally Lipschitz functional
$V:{\mathcal C}\to R^+$, the derivative of the
functional $V$, $D^+V:{\mathcal C}\times R^m\to R^{\star}$, is defined (in the Driver's form,
see Driver, 1962, Pepe, 2007 (A)), for $\phi\in {\mathcal C}$, $d\in R^m$, as
\begin{equation}\label{robderivatatotale}
 D^+ V(\phi,d)=
\limsup_{h\to 0^+} {1\over h} \left ( V(\phi_{h,d})-V(\phi) \right
),
\end{equation}
where: for $0<h<\Delta$, $\phi_{h,d} \in {\mathcal C}$ is given
by
\begin{eqnarray}\label{robfictitioussystemneutralold}  \phi_{h,d} (s) =  \left \{ \begin{array}{cc} \phi(s+h),
\qquad s\in [-\Delta,-h];     \\ {\cal
D}\phi+f(\phi,d)(s+h) -{\cal
D}\phi_{s+h}^{\star}+\phi(0),   \qquad  s\in (-h,0];
\end{array}
\right .
\end{eqnarray}
for $0<\theta\le h$, $\phi^{\star}_{\theta} \in {\mathcal C}$ is given by (\ref{phithetastarold}).

\noindent
Let us here recall that, for any locally Lipschitz functional $V:{\mathcal C}\to R^+$, the following result holds (see Pepe, 2007 (A))
\begin{equation}\label{robkrasovskiiequivdriver}
\limsup_{h\to 0^+}\frac{V(x_{t+h})-V(x_t)}{h}=D^+V(x_t,u(t)), \qquad t\in [0,b), \ a.e.,
\end{equation}
where $x_t$ is the solution of $(\ref{NFDEinput})$ in a maximal time interval $[0,b)$, $0<b\le+\infty$.

\begin{theorem}\label{GESISS} Let the system described by (\ref{NFDEinput}), with $u(t)\equiv 0$, be $0-$GES. Let the map $f$ in (\ref{NFDEinput}) satisfy the following hypotheses:
\begin{itemize}
\item[i)] there exists a positive real $L_0$ such that, for any $\phi_i\in {\mathcal C}$, $i=1,2$, the inequality holds
\begin{equation}
\vert f(\phi_1,0)-f(\phi_2,0)\vert \le L_0\Vert \phi_1-\phi_2 \Vert;
\end{equation}
\item[ii)]
there exists a function $L$ of class ${\mathcal K}$ such that, for any $\phi\in {\mathcal C}$, for any $u \in R^m$, the inequality holds
\begin{equation}\label{lipfinputu}
\vert f(\phi,u)-f(\phi,0)\vert \le L(\vert u\vert)
\end{equation}
\end{itemize}
Then, the system described by (\ref{NFDEinput}) is ISS.
\end{theorem}

\noindent

\section{Conlcusions}

In this paper we have dealt with converse Lyapunov-Krasovskii theorems for time-invariant systems described by NFDEs in Hale's form, with linear, strongly stable difference operator. We have proved that the well known Lyapunov-Krasovskii conditions, sufficient for the $0-$GAS property of these systems, are also necessary. Moreover, we have given necessary and sufficient Lyapunov-Krasovskii conditions for the $0-$GES propery. Finally, we have shown a link between the 0-GES propery and the ISS property. Future investigation will concern the case of nonlinear difference operator, namely ${\mathcal D}\phi=\phi(0)-g(\phi)$, $\phi\in {\mathcal C}$, with $g$ nonlinear, involving discrete as well as distributed time-delays (see Pepe, Karafyllis \& Jiang, 2008, Melchor-Aguilar, 2012). Notions for the nonlinear difference operator such as, for instance, $g(\phi)$ independent of $\phi(0)$, ${\phi\in \mathcal C}$ (see Definition 5.1, p. 281, in Hale \& Lunel, 1993), ISS (see Pepe, Jiang \& Fridman, 2008), incremental GAS, incremental ISS (see Angeli, 2002, as far as systems described by ODEs are concerned), may be instrumental in order to further extend the results presented here to systems described by more general NFDEs. Another topic which will be considered concerns the converse Lyapunov-Krasovskii Theorem for the ISS of systems described by NFDEs in Hale's form, for which sufficient Lyapunov-Krasovskii conditions are provided in (Pepe, 2007, (A), Pepe, Karafyllis \& Jiang, 2008). We believe that a converse Lyapunov-Krasovskii theorem for the notion of robust $0-$GAS (see Lin, Sontag \& Wang, 1996, as far as systems described by ODEs are concerned) may be instrumental for deriving the converse Lyapunov-Krasovskii theorem for the ISS of these systems, as it happens for systems described by RFDEs (see Karafyllis, Pepe \& Jiang, 2008 (A), Karafyllis, Pepe \& Jiang, 2008 (B)).

\end{document}